On the Curved Patterns in the Graphs of PPTs†

Revised

James M. Parks
Dept. of Mathematics
SUNY POTSDAM
Potsdam, NY USA
parksjm@potsdam.edu
April 19, 2021

**Abstract**

On his website *Pythagorean Right-Angled Triangles*, R. Knott gives a discrete plot of all Primitive Pythagorean Triples (PPTs) with legs up to length *10,000*, which he constructed using Mathematica. The patterns are very interesting, and suggest conic sections. We show that they indeed are on parabolic curves which follow, in a natural way, from the mathematics of the subject matter, and these curves can be used to construct Knott's graph.

## I. PRELIMINARY DEFINITIONS and REMARKS

These notes are written for beginning mathematicians of all ages. They have a few extra examples and discussions of the basics. Those bits can be skipped over by the more mature readers. I would be most interested in hearing from those readers who have thoughts about the material.

### Definitions

**1.a.** *A Pythagorean Triple (denote PT) is a set of 3 positive integers (a,b,c) which satisfy the Pythagoras equation* $a^2 + b^2 = c^2$*.*
*For example, (3,4,5) and (6,8,10) are PTs (Pythagorean Triples).*

**1.b.** *Given a PT (a,b,c), if the integers a, b, c are relatively prime, gcd{a,b,c} = 1, then we call the triple (a,b,c) a Primitive Pythagorean Triple (denote PPT).*
*For example, the PT (3,4,5) is a PPT, but the PT (6,8,10) is not.*

We are interested in analyzing a discrete graph of PPTs, which is construct as follows. Given a PT *(a,b,c),* with $a < b < c$, choose the ordered pair *(a,b)* and plot this in the *1st* quadrant of the *xy-plane,* Fig.1*. Then the right triangle *ABC* is congruent to the triangle formed by the origin *O*, the point *(a,0)*, and the point *(a,b)*, by Euclid's Side-Side-Side Proposition. Thus the point *(a,b)* determines a congruent copy of *ΔABC.* Another congruent copy of *ΔABC* is determined by the plot of the point *(b,a),* Fig.1.

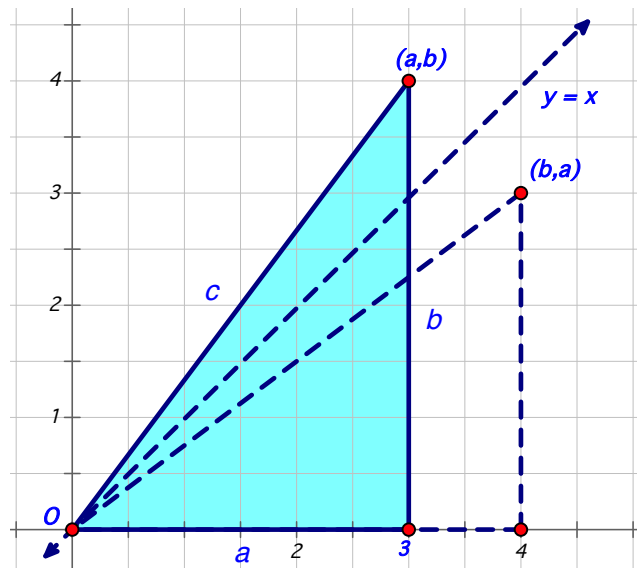

*Figure 1*

†A preliminary version of this paper appeared in [5].
*The figures in this paper were constructed using *Sketchpad* [8] and *Mathematica* [4].

## II. GRAPHS and CURVES

If all PPTs with legs *a* and *b* less than *10,000* are plotted, Fig. 2 results [3] [4] [9]. The versions [3] and [4] were constructed using *Mathematica* with different programs.

The red region is the set of points *(a,b)*, where $a < b$, and the black region is the set of points *(b,a)*. Thus the two regions are symmetric about the line $y = x$.

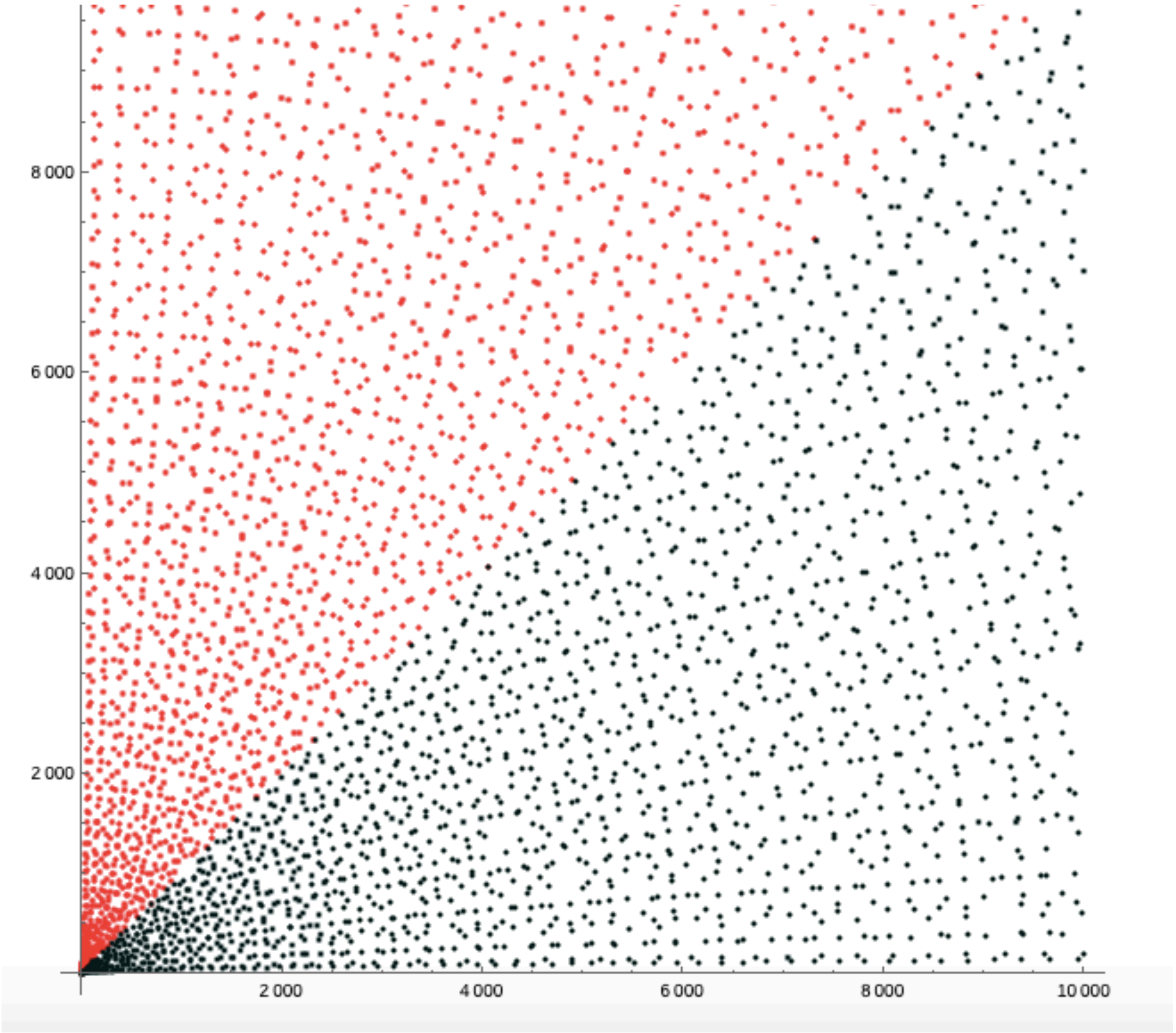

*Figure 2*

Consider the following list of the first *18* PPTs ordered by the size of the short leg *a*.

| | | | | | |
|---|---|---|---|---|---|
| ***(3,4,5)*** | ***(5,12,13)*** | ***(7,24,25)*** | *(8,15,17)* | ***(9,40,41)*** | ***(11,60,61)*** |
| *(12,35,37)* | ***(13,84,85)*** | ***(15,112,113)*** | *(16,63,65)* | ***(17,144,145)*** | ***(19,180,181)*** |
| *(20,21,29)* | *(20,99,101)* | ***(21,220,221)*** | ***(23,264,265)*** | *(24,143,145)* | ***(25,312,313)*** |

*Figure 3*

It should be obvious that if you are searching for PPTs, then the *3* terms in a triple cannot all be even, and in fact no *2* of them can be even. Likewise, the *3* terms cannot all be odd. This leaves us with the only possibility being *1* even number and *2* odd numbers. Since the legs cannot both be of odd length, one must be odd and one must be even, and the hypotenuse must therefore always be of odd length.

It can also be shown that if *(a,b,c)* is a PPT, then exactly one of *a, b* is divisible by *3*, exactly one of *a, b* is divisible by *4* (possibly the same one), and exactly one of *a, b, c* is divisible by *5*, [7].

If *a* and *c* are primes, call *(a,b,c)* a *Prime PPT.*

Thus if $a > 5$, then *b* is a multiple of *60.* For example *(11,60,61), (19,180,181)*, and *(29,420,421)*. If it is only known that *a* is a prime, this places a special restraint on the PPT, which will be discussed below.

The PPTs listed in bold print above are meant to attract your attention. These *12* particular PPTs all have 'a' terms which are odd numbers, and they all have the form *(a,b,b+1),* that is $c = b+1$.

The coordinates of these *12* are obviously relatively prime, since consecutive numbers cannot have any common divisors. But they also satisfy the equivalent equation $b = (a^2 - 1)/2$, and therefore $b+1 = (a^2 + 1)/2$. This means that $a^2$ (and thus *a*) must be an odd integer, so that *b* is then an even integer.

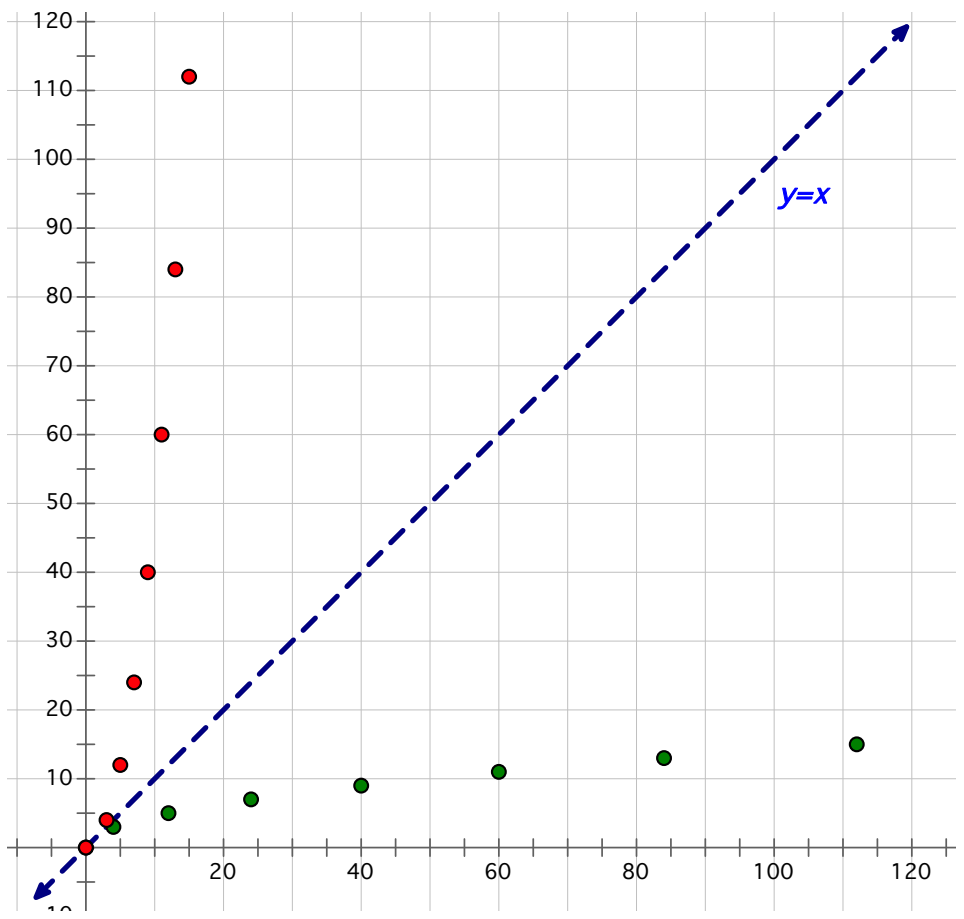

*Figure 4*

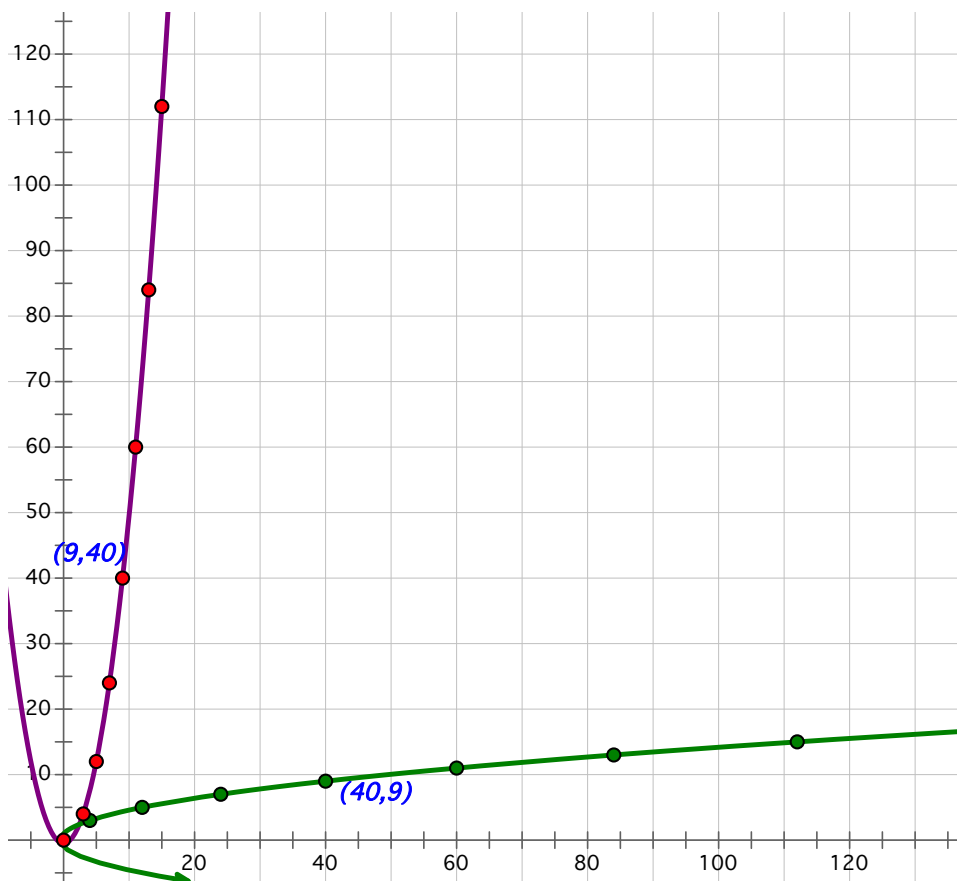

*Figure 5*

This also means that if an odd positive integer '*a*' is given, then the triple *(a,b,b+1)* is a PPT, whenever $b = (a^2 - 1)/2$. In fact these are the only PPTs which have primes in the '*a*' coordinate (more on this later). This also determines a one-to-one correspondence between the set of odd positive integers '*a*' and the set of PPTs which have the form *(a,b,b+1)*, with $b = (a^2 - 1)/2$.

This is not a new result. According to Proclus (410-485CE), these triples were known by the Pythagoreans (570-495BCE), and probably before this period [1].

A plot of these PPTs in the *xy-plane* is in Fig. 4. It suggests parabolic shaped curves opening about the positive *x,y-axes.* The equations of these parabolas: $y = (x^2 - 1)/2$, and $x = (y^2 - 1)/2$. The first parabola opens about the positive *y-axis*, with focus at the origin *O*, and vertex *(0,-1/2)*. The second parabola opens about the positive *x-axis,* has focus at *O*, and vertex *(-1/2,0),* Fig. 5.

Looking back at the list in Fig.3, notice that several triples have the form *(a,b,b+2).* Denote these PPTs by $d = 2$, since the difference $d = c - b = 2$. Similarly, we will denote the PPTs studied above, with $c = b+1$, by $d = 1$.

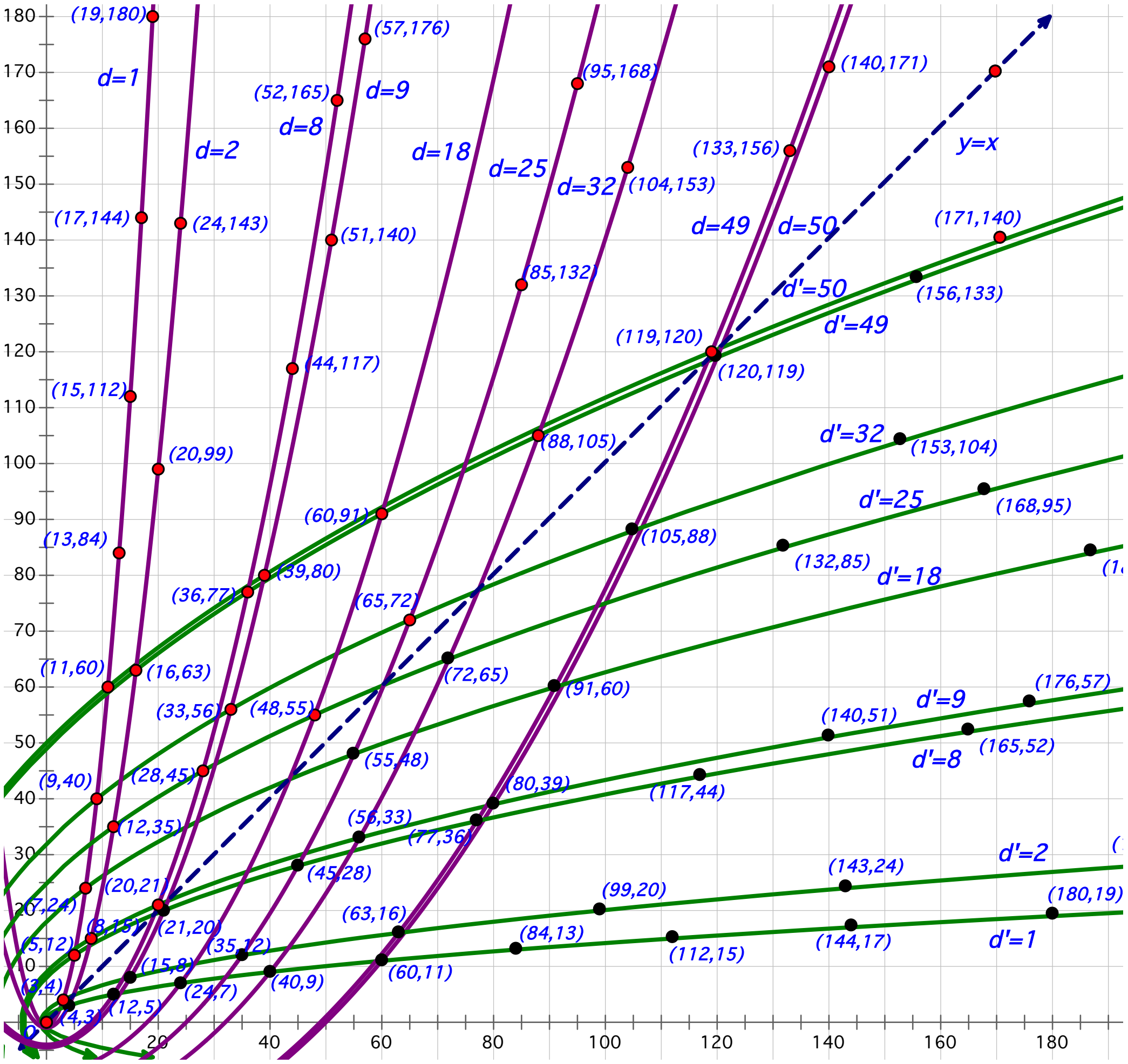

*Figure 6*

The set of PPTs which satisfies $d = 2$ must also satisfy the equivalent equation $b = (a^2 - 4)/4$. This means that $a^2$, and thus *a,* must be even, and *b* is odd. The numbers a, *b, c* are relatively prime, since the odd terms *b, c,* are *2* units apart.

Proclus attributes this result to Plato (429-347BCE) [1].

Also, from the list in Fig. 3, with $d = 2$, note that *a* is even, and $a = 0 \bmod 4$. But if *a* is even, and not a multiple of *4,* then $a = 2 \bmod 4$, and the triple is not a PPT. For if $a = 0 \bmod 4$, then $a = 4k$, $k > 0$, and $(a, b, b+2) = (4k, 4k^2 - 1, 4k^2 + 1)$, a PPT. But if $a = 2 \bmod 4$, then $a = 4k + 2$, and $b = 4k^2 + 4k$, so $b + 2 = 4k^2 + 4k + 2$, all multiples of *2*, and (a,b,c) is a PT, but not a PPT. Thus every other even term '*a*' determines a PPT. The curves which contain these PPTs are the parabolas: $y = (x^2 - 4)/4$, and $x = (y^2 - 4)/4$.

In general, the parabolas for arbitrary positive integer values of d are given by the equations $y = (x^2 - d^2)/(2d)$, and it's reflection about $y = x$, $x = (y^2 - d^2)/(2d)$, Fig. 7.

However, a lot of these values of *d* do not determine PPTs at all [7]. For example, the value $d = 3$ does not determine any PPTs, for all such PTs have a common factor of *3* in their coordinates. A similar result holds for the values *d = 4, 5, 6, 7, 10, 11, 12, 13, 14, 15, 16,* and *17*. The next PPTs occur when *d* = *8* and *9*, and the *d* = *8* case has patterns similar to those for *d* = *2.*

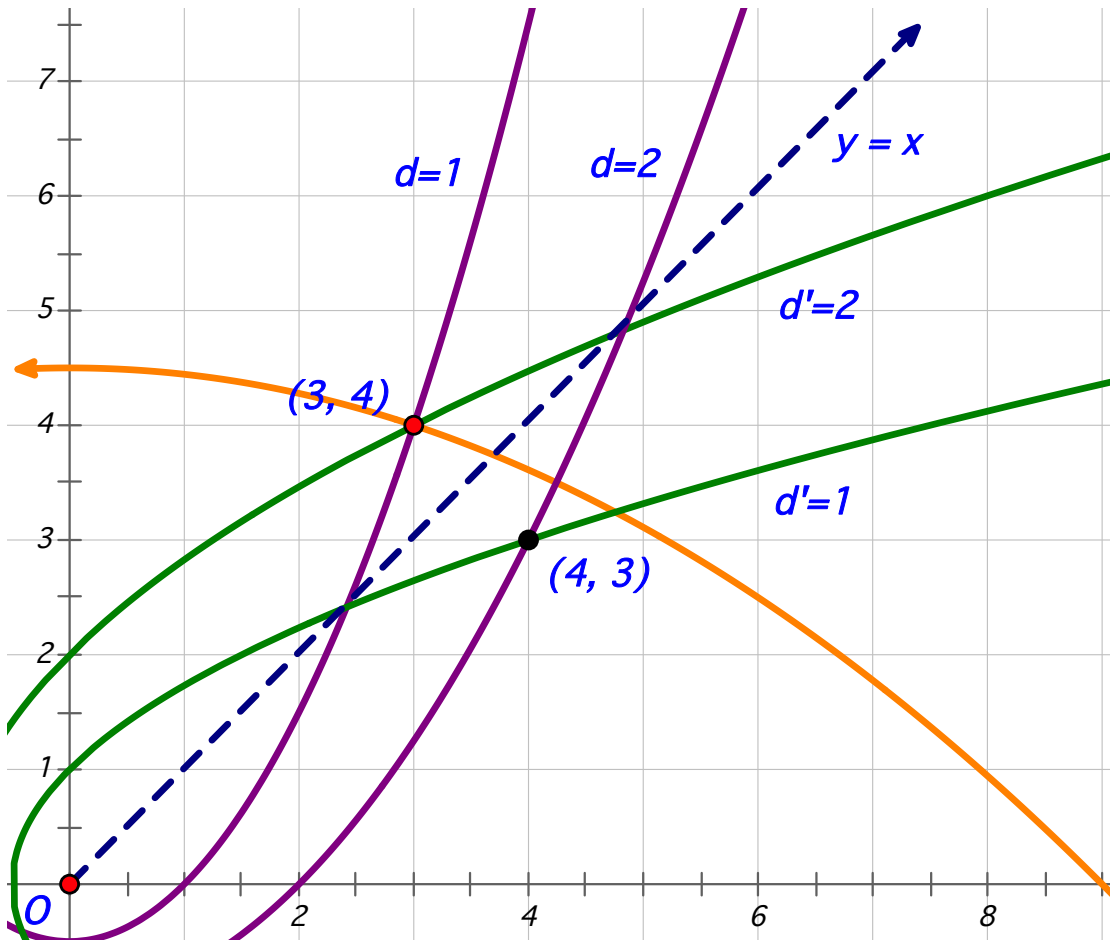

*Figure 7*

There is another *d*-value associated with the PPT *(a,b,c)* which is determined by the form *(b,a,c), a < b,* namely *d' = c - a*, where *b* is now the first coordinate of the pair.

The points *(a,b)* and *(b,a)* will occur at the intersection of the corresponding parabolas for *d* and *d'*. For example, the graph here shows the point *(3,4)* at the intersection of the *d* = *1* and *d'* = *2* parabolas, and *(4,3)* at the intersection of the *d'* = *1* and *d* = *2* parabolas, Fig.7.

All of these parabolas have their foci at the origin *O*, since the *y* coordinate of the focus point of parabolas of this type occurs at the point *x*, where $D_x(y) = \pm 1$. This happens when $x = \pm d$, and this is exactly when $y = 0$.

It was alluded to above that there is a special result if *a* is an odd prime *p, a* = *p*. For then *(p,b,b+d)* is a PPT, iff *d* = *1*, and thus $b = (p^2 - 1)/2$. This follows by the Pythagorean equation, for if *(p,b,b+d)* is a PPT, then $p^2 = 2pd + d^2 = (2p + d)d$, and $d|p^2$, thus $d|p$. This means *d* = *1*, since *d* = *p* leads to a contradiction. Hence *(p,b,b+1)* is a PPT [6].

Thus, if *a* is an odd prime *p*, then *(a,b,c)* = *(p,b,b+1)*, where $b = (p^2 - 1)/2$. But it is not a prime PPT, whenever *a* > *5*, and *b* is not divisible by *60*, by the observations on p. 3. For example, *(51,1300,1301)* is not a prime PPT, even though *51* is prime, since *1300* is not divisible by *60*.

The sequence of *d*-values which determine PPTs begins with the numbers *1, 2, 8, 9, 18, 25, 32, 49, 50, …*, see the OEIS sequence A096033 [2]. Note that the odd numbers are all squares of the consecutive odd positive integers, and the even numbers are all *2* times the square of the consecutive positive integers. We call these numbers the *allowable values of d*. The graph of a few of the representative points and parabolic curves for a few of these allowable *d* and *d'* values is shown in Figs. 7, 8.

The plot points on each parabola satisfy the spacing factor *a = r modm*, with some side conditions for certain *d* values. For example, if *d = 1, a = 1 mod 2*, and if *d = 2, a = 0 mod4.*

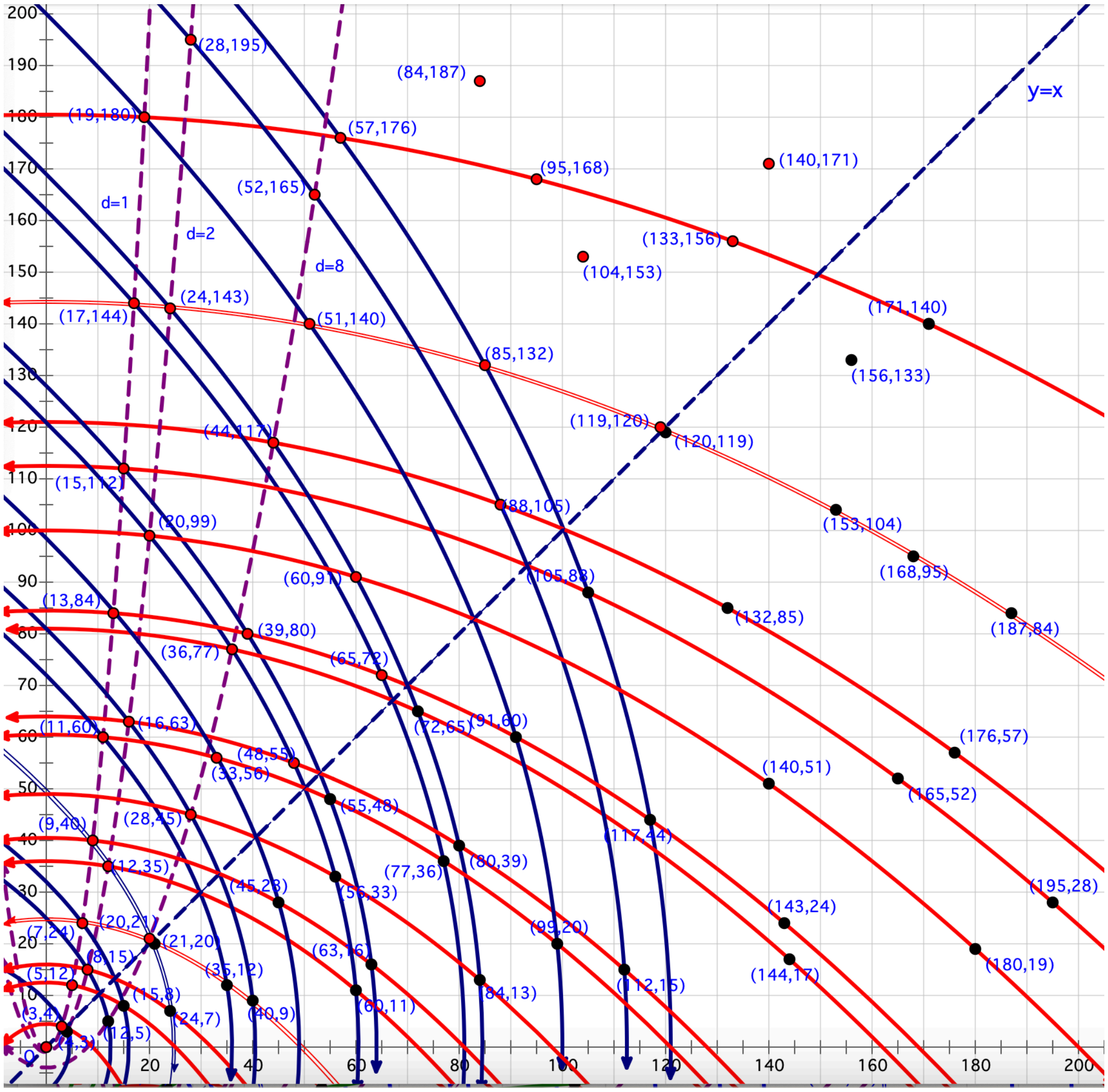

*Figure 8*

## III. PROPOSITIONS

Several of the above results are collected in the following Proposition.

**Proposition 1.** ***PPTs with Plots on Parabolas Opening on Positive x,y-axes.***
*If (a,b,c) is a PPT with a < b < c, then (a,b) is at the intersection of graph of the parabolas with equations* $y = (x^2 - d^2)/(2d)$, *and* $x = (y^2 - d'^2)/(2d')$, *for the allowable values* $d = c - b$, *and* $d' = c - a$.

It also follows that *(b,a)* is at the intersection of the reflection of these parabolas about $y = x$, which have equations $x = (y^2 - d'^2)/(2d')$, and $y = (x^2 - d^2)/(2d)$.

On the Main Graph, Fig. 3, notice that there also appear to be parabolic shaped curved patterns which open about the negative *x, y-axes*. These curves are mentioned in [3]. The equations for these curves are determined in a different way from those above.

First, consider the parabolic curves which open about the negative *y-axis*. The general form for these equations is $y = Ax^2 + B$, where $A < 0$. If *2* points on the curve are known, then the equation can be determined. For example, a possible curve might exist which contains the points *(13,84), (39,80),* and *(65,72).* The equation of this parabola is therefore $y = - x^2/(2x13^2) + 13^2/2$. One of the pairs in this example, (39,80), is on the $d = 9$ parabola, since the PPT is *(39,80,89).* Here $\sqrt{d} = 3$, and $39 = 3x13$. Similarly for the pair (*65,72)*, which is on the $d = 25$ parabola, and $65 = 13x5$, Fig. 9.

So, if the PPT plot point *(a,b)* has coordinate $a = 13m$, and the point is on some *d* parabola, *d* odd, then $m = \sqrt{d} = s$.

In general, if *a* is odd, and *(a,b)* is on some *d* parabola, *d* odd, then $a = a_1s$, where $a_1$ is the first coordinate of the corresponding PPT pair plotted on the $d = 1$ parabola (recall coordinate *a* cannot be prime if $d > 1$).

Thus the general form for the equations of these parabolas which open about the negative *y-axis*, is $y = - x^2/(2a_1^2) + a_1^2/2$, where $a_1$ is the value of the *x-coordinate*, on the corresponding parabola $d = 1$, (on the LH red dashed curve in the graph just above), and $a_1$ divides all of the first coordinates of the PPT plots on that curve.

Similarly to the case of the parabolas which open up about the positive *y-axis*, all of these parabolas also have their foci at the origin *O*. This is because the *y-coordinate* of the focus point of parabolas of this type occurs at the point *x,* where $D_x(y) = +1$, and this happens when $x = a_1^2$, which is exactly when $y = 0$.

In the graph in Fig. 8, notice that some sequences of points start with coordinate a even. For example, the pair *(20,99)* is on the $d = 2$ parabola, and is the plot of the PPT *(20,99,101).* The equation for the parabola which contains the plot of this point is $y = -x^2/20^2 + 20^2/4$. For another example with coordinate a even, consider the pair *(28,45).* This pair is on the $d = 8$ parabola, and is from the PPT *(28,45,53).* The equation for this parabola is $y = - x^2/(28^2/4) + 28^2/16$.

Thus the general form for the equations of these parabolas, where *a* is even, and they open about the negative *y-axis*, is $y = - x^2/(2a^2/d) + a^2/(2d)$. This is given that the series of pairs has even 1st coordinates, all multiples of *a*, where *(a, b)* is on either the $d = 2$ or $d = 8$ parabola. The focus point of these parabolas is the origin *O*, by a similar argument as given above for the case $d = 1$.

Every pair *(a,b)* on the $d = 2$ and $d = 8$ parabolas determines a parabola in this set. Thus every plot point from a PPT on the $d = 1$, $d = 2$, and $d = 8$ parabolas determines a parabola with equation $y = - x^2/(2a^2/d) + a^2/(2d)$.

These parabolic curves are shown in orange in the graph in Fig. 9. When these sets of parabolas are reflected about the line $y = x$, the corresponding parabolas are for the *d'-values*. These parabolas open about the negative *x-axis*, and are shown in dark blue in the graph, Fig. 9.

On the other hand, every parabola shown in the graph above, which opens about the negative *y-axis*, has a point *(a,b)* from one of the parabolas for *d = 1, 2,* or *8*. This follows from the observation that these parabolas are the closest ones to the *y-axis* for even or odd values of *d*.

Given a point *(a,b)* on the discrete graph Fig.2, it was shown above that one parabola which contains *(a,b)* has equation $y = (x^2 - d^2)/(2d)$. If *a* is odd, then either *(a,b)* is on the $d = 1$ parabola, $y = (x^2 - 1)/2$, and $b = (a^2 - 1)/2$, or *(a,b)* is on a *d*-parabola with $d > 1$, $d$ odd, and $b < (a^2 - 1)/2$. Plus, $a = ta_1$, and $d = t^2$, $t$ odd, where $a_1$ is the smallest odd number such that $a_1|a$, and $b < b_1 = (a_1^2 - 1)/2$.

Thus $(a_1,b_1)$ is on the $d = 1$ parabola, and $(a,b)$ is on the parabola which opens down about the negative $x$-axis, with equation $y = - x^2/(2a_1^2) + a_1^2/2$.

For example, $(105,88)$ is a point on the graph of all PPTs, with $c = 137$, so $d = 49$, and $d' = 32$. Thus $(105,88)$ is on the $d = 49$ parabola $y = (x^2 - 49^2)/98$. Since $b_1 = (a_1^2 -1)/2 > 88$, where $a_1$ divides $105$, solving this inequality for $a_1$ determines that $a_1 > \sqrt{177} \sim 13.3$. Since $105 = 3x5x7$, it is suspected that $(a_1,b_1) = (15,112)$. This is indeed true, as $c_1 = 113$. Thus $(105,88)$ is also on the parabola with the equation $y = -x^2/(2\cdot 15^2) + 15^2/2 = -x^2/450 + 225/2$, as can be easily checked, Fig.10.

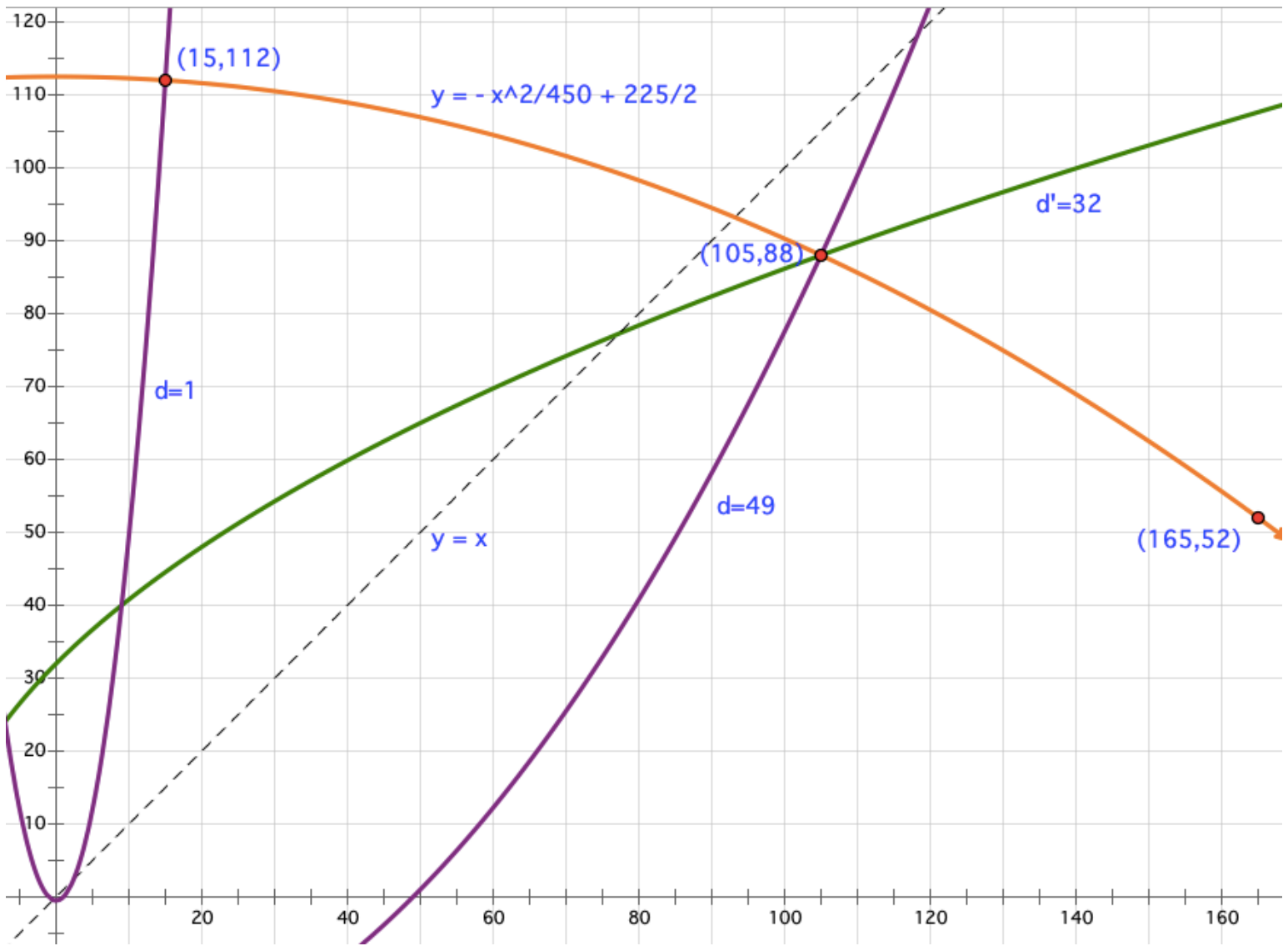

*Figure 9*

If $a$ is even, then a similar argument applies, where $b_1 = (a_1^2 - 4)/4$, or $(a_1^2 - 64)/16$, depending on whether $d = 2$ or $8$, respectively.

For example, let $(a,b) = (140,51)$. Then $4|140$, but $8$ does not. Also $c = 149$, so $d = 98$, and $(a,b)$ is on the parabola with equation $y = x^2/196 - 49$. Since $140 = 2^2x5x7$, it follows that $b_1 = (a_1^2 - 4)/4 > 51$ iff $a_1 > \sqrt{208} \sim 14.4$. Let $a_1 = 20$, the smallest divisor of $140$ greater than $14.4$, which is also divisible by $4$. Then $b_1 = 99$, $c_1 = 101$, $d = 2$, and the point $(20,99)$ is on the parabola with equation $y = - x^2/400 + 100$. A similar argument holds if $8|a$.

Given a point $(b,a)$ on the graph of all PPTs, with $a < b$, $a$ odd, this point is on the graph of a parabola about the negative $x$-axis, with equation $x = - y^2/(2a_1^2) + a_1^2/2$, $a_1$ as above, and $d' = 1$. For the reflected point $(a,b)$ is on the graph of the reflected parabola about the negative $y$-axis, $y = - x^2/(2a_1^2) + a_1^2/2$, with $d = 1$. The techniques above, and a reflection of this parabola about $y = x$ gives the desired equation. Thus these results determine the following Proposition.

**Proposition 2.** ***Parabolas About the Negative y-axis Determined by a Single Point***

*Given a PPT plot point $(a,b)$, a parabola with equation $y = - x^2/(2a_1^2/d) + a_1^2/(2d)$ is determined which contains it, and which opens about the negative y-axis. When $a$ is odd, and $a_1$ is the smallest odd number such that $a_1 \mid a$, and $b < b_1 = (a_1^2 - 1)/2$, then $(a_1,b_1)$ is on the $d = 1$ parabola. For $a$ even, $(a_1,b_1)$ is on the $d = 2$ parabola, if $a_1$ is even, and $b_1 = (a_1^2 - 4)/4$, or on the $d = 8$ parabola, if $a_1$ is even, and $b_1 = (a_1^2 - 64)/16$.*

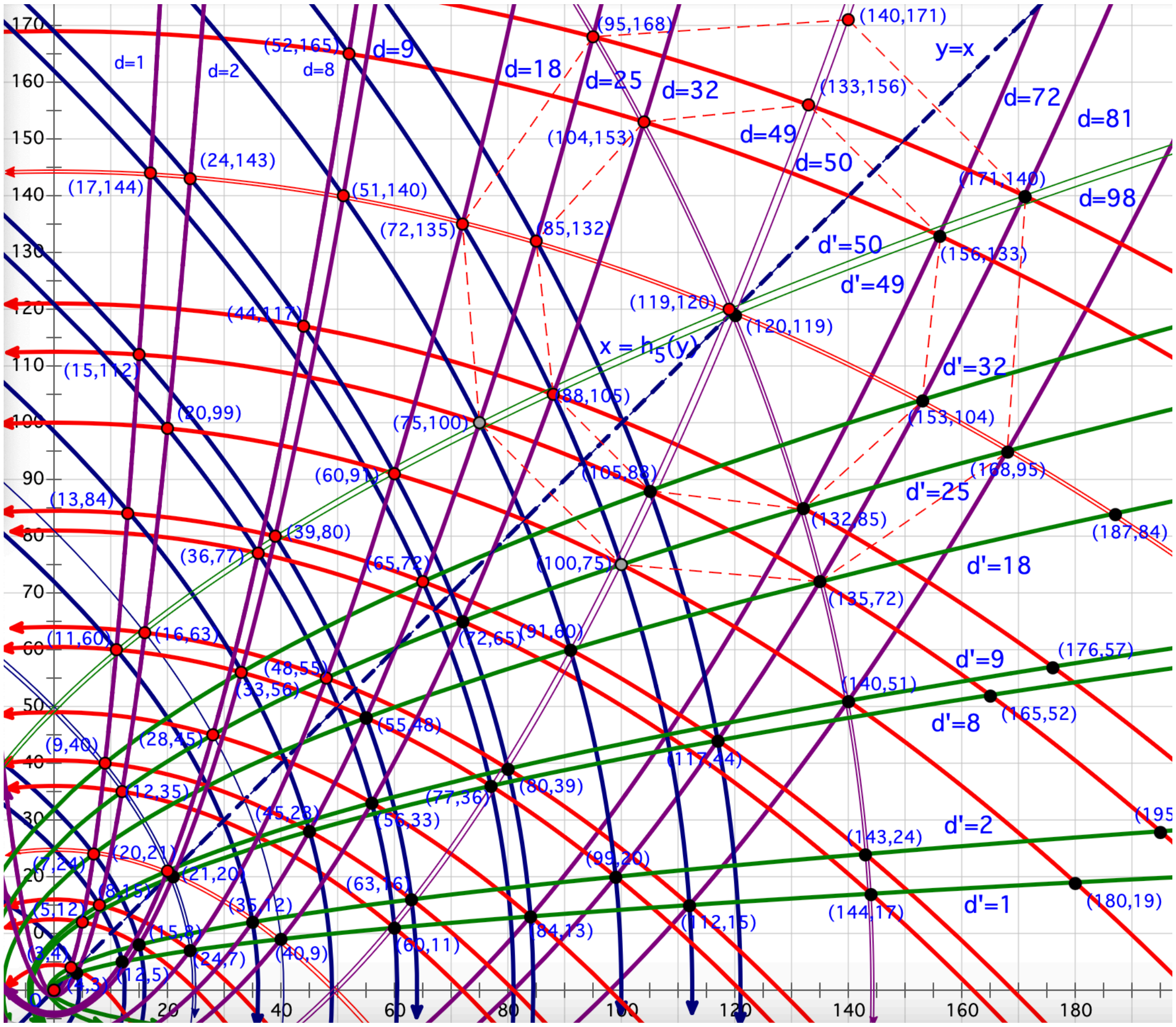

*Figure 10*

In addition to the above, for the parabolas with allowable d values, which open about the positive y-axis, they intersect the parabolas which open about the negative y-axis at right angles at the PPT plot points of intersection. This is easily seen by calculating the slopes of the tangent lines to the parabolic curves at the point of intersection. For example, since the $d = 1$ parabola has equation $y_1 = (x^2 - 1)/2$, let $(a_1, (a_1^2 - 1)/2)$ be the point of intersection of this curve with the parabola $y_2 = - x^2/(2 \cdot a_1^2) + a_1^2/2$. Then, using the derivatives of the parabolas, the slope of the tangent line to parabola $y_1$ at the point $(a_1, (a_1^2 - 1)/2)$, is $m_1 = D_x y_1 = a_1$, and the slope of the tangent line to parabola $y_2$ at this same point is $m_2 = D_x y_2 = -1/a_1$. Hence, $m_1 m_2 = -1$, and the parabolas are perpendicular at $(a_1, (a_1^2 - 1)/2)$.

In general, let $d$ and $a$ be odd positive integers, $d$ allowable, and let $(a,b)$ be a point on the $d$ parabola, which has equation $y_1 = (x^2 - d^2)/(2d)$. Then $(a,b)$ is on a parabola which opens down about the negative $y$-axis, with equation $y_2 = - x^2/(2 \cdot a_1^2) + a_1^2/2$, where $(a_1,b_1)$ is on the, say $d_0 = 1$ parabola, $a = sa_1$, and $d = s^2$.

It is claimed that $y_1$ and $y_2$ are perpendicular. Observe, the slopes satisfy $m_1 = D_x y_1 = a/d = a_1/s$, and $m_2 = D_x y_2 = - a/a_1^2 = - s/a_1$, at $x = a$, so $m_1 m_2 = - 1$. Thus the parabolas are perpendicular at $(a,b)$. Similar arguments apply for $d$ and $a$ even.

In the cases for pairs of parabolas which open about the positive *x-axis*, and the negative *x-axis*, these parabolas are perpendicular, since they are reflections about $y = x$ of the previous case of parabolas which opened about the *y-axis*.

These results are summarized in the following Proposition.

**Proposition 3.** ***Perpendicular Parabolas on One Axis***

*If (a,b) is a point on the graph of all PPTs, which is at the intersection of a pair of parabolas with allowable d-values, one which opens up and one which opens down about the y-axis, or one which open to the right and one which opens to the left about the x-axis, then they meet at right angles at the point (a,b).*

**IV. The OEIS SEQUENCE A096033**

Now, why exactly are the allowable values of d equal to those numbers given in the OEIS sequence A096033 [6]? As noted above, the sequence is made up of the squares of the odd positive integers: *1, 9, 25, 49, 81, …* , and those even positive integers which are twice the squares of the consecutive positive integers: *2, 8, 18, 32, 50, 72, …* , and so on.

First, suppose $d$ is odd, and $d \neq s^2$, $s$ odd. Note that, since $a^2 = c^2 - b^2 = (c - b)(c + b) = d(c + b)$, then $d|a$ for a given associated PT $(a,b,c)$, with $d = c - b$. Assume $a$ is an odd integer, then empirically it follows that $a = d \bmod 2d$, for $a$ sufficiently large, so $a = (2k +1)d$, $k \geq 1$. This means $b = (4d^2(k^2 + k) + d^2 - d^2)/(2d) = 2(k^2 + k)d$, and $c = (2(k^2 + k) +1)d$. But then all coordinates have the common factor $d$, so $(a,b,c)$ is not a PPT, and thus $d$ is not on the list of allowable values A096033.

On the other hand, if $d$ is odd, and $d = s^2$, $s$ odd, then for a corresponding PT $(a,b,c)$, $a = s \bmod 2s$, by empirical data. Thus $a = (2k + 1)s = ms$, an odd integer, where $k \geq 1$. Since $a^2 = c^2 - b^2 = (c - b)(c + b)$, assume $m > s$, and let $m^2 = (c + b)$ and $s^2 = (c - b)$. Then $b = (m^2 - s^2)/2$, an even integer, and $c = (m^2 + s^2)/2$, an odd integer. It is claimed that there are no common factors here, since m and s do not divide either $m^2 - s^2 = (m - s)(m + s)$, an even integer, or $m^2 + s^2$, an even integer which satisfies $m^2 < m^2 + s^2 < 2m^2$. So $(a,b,c)$ is a PPT, and $d$ is in A096033.

Suppose $d$ is even, and $d \neq 2s^2$, $s$ a positive integer, then $d = 2k$, for some $k$, a positive integer, and $k \neq s^2$. If $(a,b,c)$ is a corresponding PT, then $a^2 = 4k(b + k)$, so $a$ is even, and $a = 2nd = 4nk$, for $n = (b + k)$. Thus $b = (4n^2 - 1)k$, and $c = (4n^2 + 1)k$. Since all terms share the common factor $k$, $(a,b,c)$ is not a PPT, and $d$ is not in A096033.

If $d$ is even, and $d = 2s^2$, there are two cases to consider:

<u>First case</u>: Let $s$ be odd, and assume for an associated PT $(a,b,c)$ that $a = 0 \bmod 4s$, by empirical data. Then $a$ is even, $a = 4sk$, for $k$ an integer, and assume $k \neq 0 \bmod s$, and $k > s$. So $b = 4k^2 - s^2$, and $c = 4k^2 + s^2$. Now $k$ does not divide $b$, since $b = (4k^2 - s^2) = 11\ (2k - s)(2k + s)$. Similarly, $s$ does not divide $b$, since $s$ does not divide $k$. Thus $(a,b,c)$ is a PPT, and $d$ is on the list A096033.

<u>Second case</u>: Let $s$ be even, $s > 4$, then, for an associated PT $(a,b,c)$, where $a = 0 \bmod 2s$ by empirical data, $a = 2sk$, for $k$ odd. Assume $k > s$, so $a$ is an even integer. Then $s$ does not divide either $b = (k^2 - s^2)$ or $c = (k^2 + s^2)$, as they are both odd. Since $k > s$, and $b = k^2 - s^2 = (k - s)(k + s)$, $k$ does not divide $b$. Thus $(a,b,c)$ is a PPT, and $d$ is on A096033.

## V. CIRCULAR PATTERNS

The many circular patterns which appear in the Main Graph, Fig.2, are also of interest. Several examples are given here which show that actual circular shapes are impossible, but the shapes are still mystifying. Most of the circular patterns are made up of 8 plot points of PPTs. Some of the points of a circular shape may be missing if the corresponding triple is a PT, but not a PTT.

In Fig. 11 notice that *2* circular patterns in a starburst pattern are visible with centers near the point *(120,120).* The octagon shapes shown highlight their graphs. The smaller inner circular pattern is made up of the plot points listed here, arranged in a ccw direction with the first *4* points corresponding to the plot points below the line $y = x$ (black), starting with *(105,88),* and the last *4* points corresponding to the reflected plot points (red), of the first *4*, ending with *(88,105).*

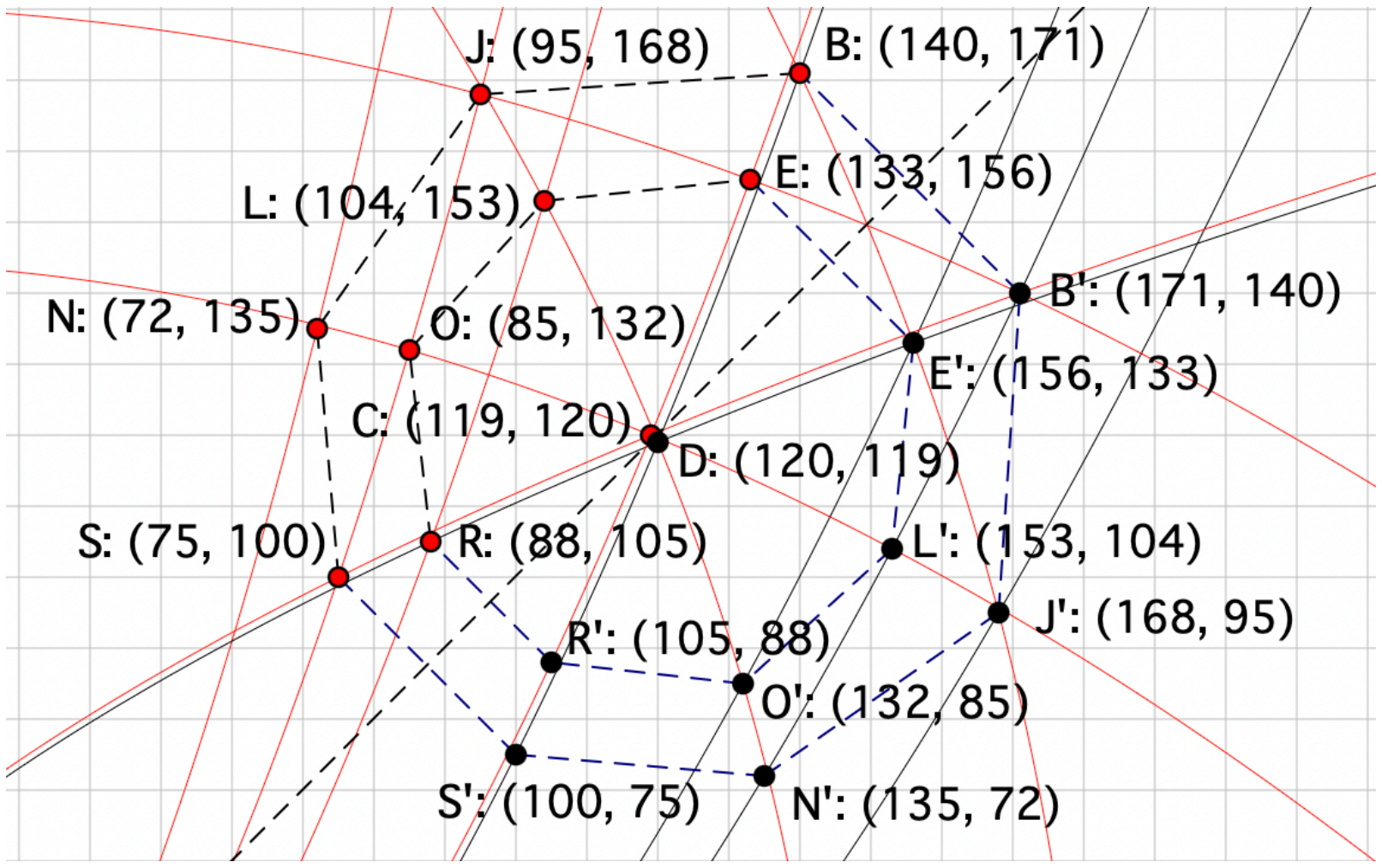

*Figure 11*

So as you travel around this circular shape starting at *(105,88)*, and ending at *(88,105),* you determine the following list of points:

***(105,88), (132,85), (153,104), (156,133)****, (133,156), (104,153), (85,132), (88,105).*

The largest prime divisors of each coordinate of each point determines the following sets of pairs, presented here in the same corresponding order as the list of point plots.

***{7,11}, {11,17}, {17,13}, {13,19},*** *{19,13}, {13,17}, {17,11}, {11,7}.*

This determines a set of links between the pairs of prime divisors of the PPTs whenever the *2nd* coordinate of a pair equals the *1st* coordinate of a following pair.

At first glance all of the links between these pair endpoints is surprising. But because of reflections about $y = x$, only *3* links are actually important, shown in the list above in bold, and equivalently in non-bold, all of the rest follow by reflections about $y = x.$

But is this result a general truth? We need more examples.

Here is another example of a circular pattern on the line $y = x,$ which has center near *(400,400),* Fig. 12:

***(377,336), (420,341), (465,368), (460,429)****, (429,460), (368,465), (341,420), (336,377).*

The corresponding largest prime devisor pairs are listed here:

***{29,7}, {7,31}, {31,23}, {23,13},*** *{13,23}, {23,31}, {31,7}, {7,29}.*

This is the same relationship that we saw in the example above, all pairs are linked.

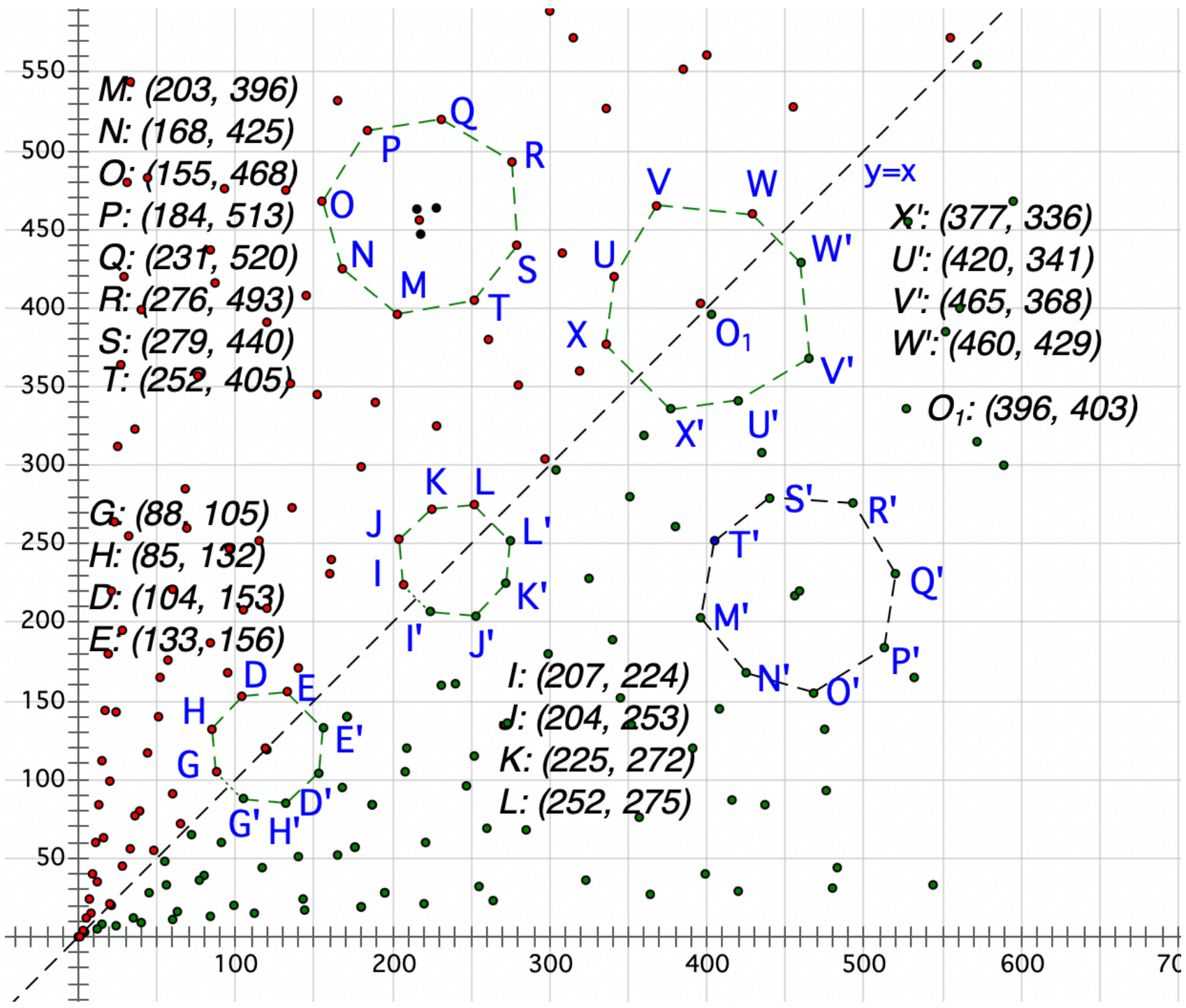

*Figure 12*

The other circular region in Fig. 12 has the following plot points:

***(100,75), (135,72), (168,95), (171,140),*** *(140,171), (95,168), (72,135), (75,100).*

The largest prime divisors of each coordinate determines the following list of pairs:

***{5,5}, {5,3}, {7,19}, {19,7}***, *{7,19}, {19,7}, {3,5}, {5,5}.*

We now see that the *2nd* pair has *2nd* coordinate *3*, and the *3rd* pair has *1st* coordinate *7*, so these two pairs of divisors and their reflected pairs are not linked.

However, another interesting connection occurs here. Since *72 = 24x3,* and *168 = 24x7,* there is an alternate link between *72* and *168,* namely *24*, the cofactor of *3* and *7.*

So we denote this special new type of link by listing *3(24)* and *7(24)* in place of *3* and *7*:

***{5,5}, {5,3(24)}, {7(24),19}, {19,7},*** *{7,19}, {19,7(24)}, {3(24),5}, {5,5}.*

Although these two plot point pairs have no prime links, they are connected, in a sense, by the cofactor *24*, which makes a new type of link.

What about examples not centered near $y = x$?

In Fig. 13, on the left, there is a circular pattern of *8* plot points centered near the point *(2200,3900)* with line segments joining them. The *5* parabolas which contain these *8* plot points satisfy the following:

| *d* | Equation | PT plot point(s) |
|---|---|---|
| *512* | $y = x^2/1024 - 256$ | *(2080,3969)* |
| *529* | $y = x^2/1058 - 264.5$ | *(2093,3876), (2139,4060)* |
| *578* | $y = x^2/1156 - 289$ | *(2176,3807), (2244,4067)* |
| *625* | $y = x^2/1250 - 312.5$ | *(2275,3828), (2325,4012)* |
| *648* | $y = x^2/1296 - 324$ | *(2340,3901)* |

The sequential list of plot points is:

*(2176,3807), (2275,3828), (2340,3901), (2325,4012),*
*(2244,4067), (2139,4060), (2080,3969), (2093,3876).*

Clearly these *8* points are not on a perfect circle, see for example the closeup view on the right-hand graph, Fig. 13.

The largest prime divisors of each point coordinate determines the following list of pairs, presented here in the same corresponding pattern as the PPTs:

*{17,47}, {13,29}, {13,83}, {31,59}, {17,83}, {31,29}, {13,7}, {23,19}.*

This list contains no links or cofactor links, which suggests that being away from the line *y* = *x* determines a very different relationship for the circular pattern of points.

A similar situation occurs with the points in the circular pattern *M', N', O', P', Q', R', S',T'* in Fig. 12.

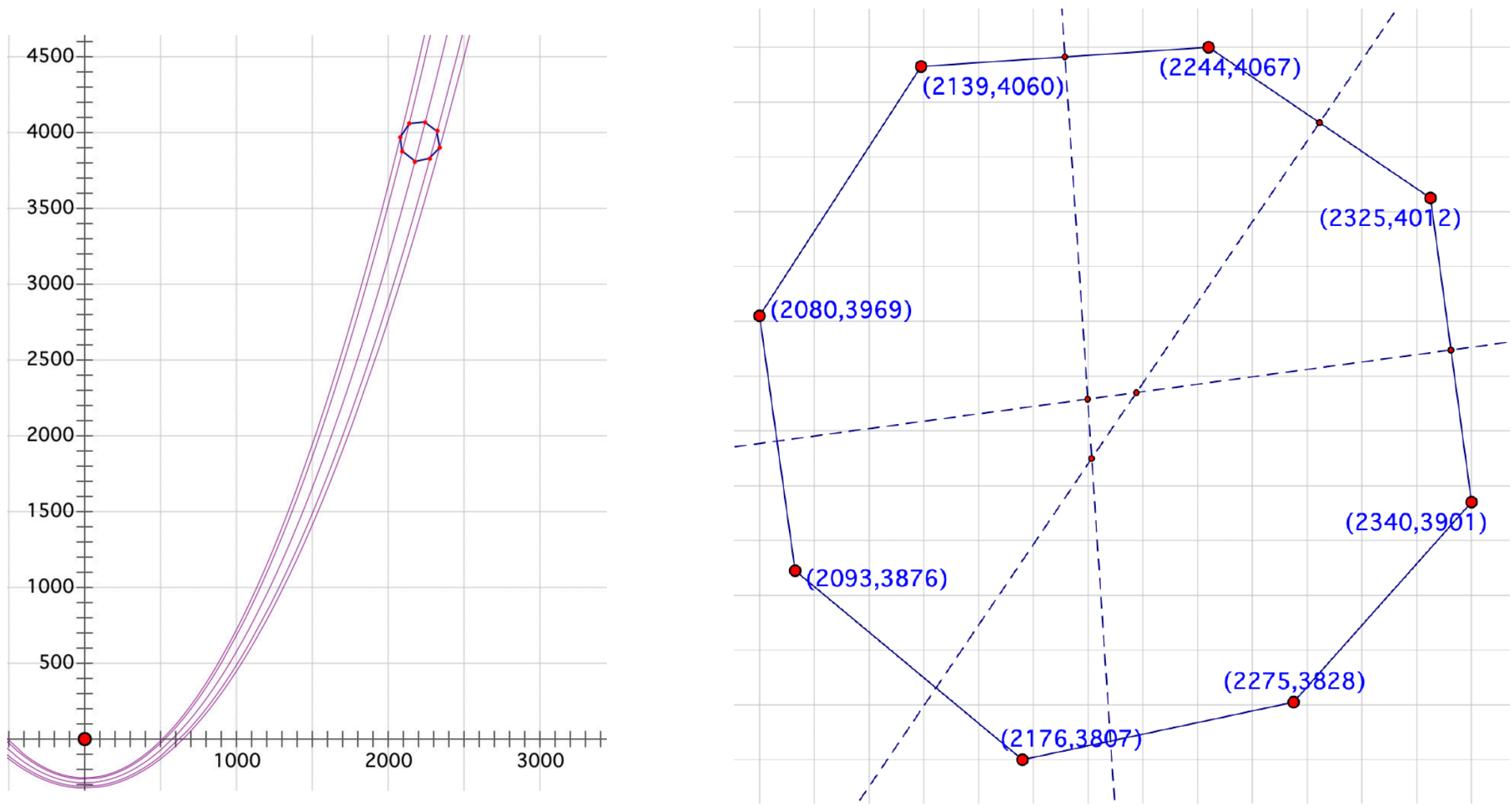

*Figure 13*

Each of the circular patterns in the Main Graph, Fig.2, appear to be a plot on at least *2* intersecting sets of *5* parabolic curves. Since these parabolas are not parallel to each other, it should not be expected that the plot points would form a perfect circle.

To see just how complicated the conditions are which the plot points in a circular pattern must satisfy, look at the set of *8* points in the circular pattern of the example near *(120,120)* in more detail, Fig.11. They must satisfy the following *5* conditions simultaneously.

- The LH plot point *(85,132)* is on the *d* = *25* parabola, where the points must satisfy: *a* = *5 mod10*, *a* not divisible by *25*.
- The next parabola is *d* = *32,* and it contains the points *(88,105)* and *(104,153)*, which both must satisfy *a* = *0 mod8, a* not divisible by *16.*
- That parabola is followed by the *d* = *49* parabola which contains the points *(105,88)* and *(133,156).* These points must satisfy *a* = *7 mod14,* that is *a* is an odd multiple of $7$, and is not divisible by $49$. This parabola also contains the central point $(119,120)$. Notice that the difference between the consecutive '$a$' points here is $14$.

- The next parabola is the $d = 72$ parabola which contains the points $(132,85)$ and $(156,133)$, which must satisfy $a = 0 \bmod 12$, that is $a$ is a multiple of $12$, and not divisible by $36$.
- The last parabola on the RH side is the $d = 81$ parabola, which contain the point $(153,104)$, and which satisfies $a = 9 \bmod 18$, that is $a$ is an odd multiple of $9$, with the except of the multiples $k$, when $k = 1 \bmod 3$.

In general, if $(a,b)$ is a point on the $d$-parabola, $d = s^2$, $s$ odd, then $a$ is an odd multiple of $s$, $a = s \bmod 2s$, and $a$ is not divisible by $d$, plus possibly some other restrictions. But, if $d = 2s^2$, $s$ even, then $a$ is a multiple of $2s$, $a = 0 \bmod 2s$, and $a$ is not divisible by $d = s^2$.

Let $(a_1,b_1)$ and $(a_2,b_2)$ be points on a $d$-parabola, then the $b$ coordinates of these points must also satisfy $b_1 = (a_1^2 - d^2)/(2d)$ and $b_2 = (a_2^2 - d^2)/(2d)$. Let $a_1$ and $a_2$ be consecutive terms on this $d$-parabola, then $a_2 = a_1 + m$. If $d = s^2$, $d$ and $s$ odd, then $a_1 = s \bmod 2s$, $a_1$ is an odd multiple of $s$, and not divisible by $d = s^2$. If $d$ is even, $d = 2s^2$, then $a_1 = 0 \bmod 2s$, so $a_1$ is a multiple of $2s$, and it is not divisible by $s^2$. Thus $b_2 - b_1 = (2a_1m + m^2)/(2d)$.

Unlike the $a$ coordinates, the $b$ coordinates are increasing in size spacing between consecutive terms, because of the '$a_1$' term in the expression for $b_2 - b_1$. Thus, two points $(a_1,b_1)$ and $(a_2,b_2)$ are near each other, if the difference $b_2 - b_1 = (2a_1m + m^2)/(2d)$ is not too large, with respect to the difference $a_2 - a_1$, relatively speaking.

For example, consider the plot points $(88,105)$ and $(104,153)$ in the example above, Fig.11. These points are on the $d = 32$ parabola, and have $a_2 - a_1 = 16$, where $a = 0 \bmod 8$, $a$ not divisible by $d$, so $b_2 - b_1 = 48$, obviously. Therefore these two points are fairly close, graphically, plus the observation that the diameter of the circular region is approximately $75$ units.

## References

[1] T. Heath, *Euclid, The Thirteen Books of the Elements,* Vol. 1, 2nd ed., Dover, 1956.

[2] *The OEIS Foundation, Inc.*, Entry A096033. Published electronically at https://oeis.org, author Lekraj Beedassy, 2004.

[3] R. Knott, *Pythagorean Right Triangles,* Section 4.4.1, http://www.maths.surrey.ac.uk/hosted-sites/R.Knott/Pythag/pythag.html#section3.4

[4] J. Parks, *Mathematica Notebook Archive*, https://www.notebookarchive.org/ppt-graph-nb--2021-06-6y23baq/

[5] _______, *Curved Patterns in the Graphs of PPTs,* Proc., 26th ATCM, Asian Technology Conference in Mathematics (2021), pp.277-286.

[6] _______, *Computing Pythagorean Triples,* arXiv:2107.06891, July 14, 2021, https://arxiv.org/pdf/2107.06891.pdf

[7] W. Sierpinski, *Pythagorean Triangles*, Dover, 2003; Republication of the original work published by the Graduate School of Science, Yeshiva U., NY, 1962.

[8] *Sketchpad* v.5.11BETA(3060), https://gsptest.scratchconsortium.com/download.html

[9] WOLFRAM Demonstrations Project, *Primitive Pythagorean Triples 1: Scatter Plot,* https://demonstrations.wolfram.com/PrimitivePythagoreanTriples1ScatterPlot/